\documentclass[12pt]{article}
\usepackage{amsthm,amsfonts,amssymb,amscd}

\begin{document}

\begin{center}
\Large \bf Birationally rigid Fano varieties
\end{center}
\vspace{1cm}

\centerline{\large \bf A.V.Pukhlikov\footnote{The author was financially supported by Russian 
Foundation of Basic Research, grants 02-15-99258 and
02-01-00441, by a Science Support Foundation grant for young researchers and by INTAS-OPEN, grant 2000-269.}} 
\vspace{1cm}

\begin{center}
Max-Planck-Institut f\" ur Mathematik \\
Vivatsgasse 7 \\
53111 Bonn \\
GERMANY \\
e-mail: {\it pukh@mpim-bonn.mpg.de}
\end{center}
\vspace{0.5cm}

\begin{center}
Steklov Institute of Mathematics \\
Gubkina 8 \\
117966 Moscow \\
RUSSIA \\
e-mail: {\it pukh@mi.ras.ru}
\end{center}
\vspace{1cm}

\centerline{February 3, 2003}\vspace{1cm}

\parshape=1
3cm 10cm \noindent {\small \quad\quad\quad \quad\quad\quad\quad
\quad\quad\quad {\bf Abstract}\newline
We give a brief survey of the concept of birational
rigidity, from its origins in the two-dimensional
birational geometry, to its current state. The main
ingredients of the method of maximal singularities
are discussed. The principal results of the theory
of birational rigidity of higher-dimensional
Fano varieties and fibrations are given and certain
natural conjectures are formulated.} 

\newpage

{

CONTENTS 
\vspace{0.7cm}

\noindent
0. Introduction
\vspace{0.3cm}

\noindent
1. The Noether theorem
\vspace{0.3cm}

\noindent
2. Fano's work
\vspace{0.3cm}

\noindent
3. The theorem of V.A.Iskovskikh and Yu.I.Manin 
\vspace{0.3cm}

\noindent
4. The method of maximal singularities
\vspace{0.3cm}

\noindent
5. Birationally rigid varieties
\vspace{0.3cm}

\noindent
6. Singular Fano varieties
\vspace{0.3cm}

\noindent
7. The relative version
\vspace{0.3cm}

\noindent
References}
\vspace{0.5cm}

\newpage

{\large\bf 0. Introduction}
\vspace{0.5cm}

\noindent
This paper is based on the talk given by the author at the 
Fano conference in Turin. The aim of the paper is to give a
brief survey of the concept of birational rigidity which
nowadays is getting the status of one of the crucial concepts
of higher-dimensional birational geometry. The Fano
conference both by definition and its actual realization
had a natural historical aspect. Therefore it seems most
appropriate to review the story of birational rigidity,
presenting the principal events in their real succession,
from the first cautious steps made in XIX century to the
modern rapid development. In this story Gino Fano himself
played a prominent part: birational rigidity was one of his
most important foresights.

For about fifty years Fano was the only mathematician in
the world engaged in the field that in his time was a real
terra incognita, three-dimensional birational geometry of
algebraic varieties which are now called Fano varieties. He
had a program of his own and he did his best to realize it,
see [F1-F3].
On this way he discovered a lot of fascinating geometric
constructions, found new approaches to investigating
extremely deep and challenging problems, made a terrific
amount of computations and guessed certain fundamental
facts. He never completed his program. But even realization
of a part of it took about thirty years of hard work of his
successors, equipped with incomparably stronger
techniques.

The author is grateful to the Organizing Committee of the
conference, namely, to Prof. Alberto Conte, Prof. Alberto
Collino and Dr. Marina Marchisio for the invitation to
give this talk and for making such a wonderful conference.

The talk was prepared during the author's stay at the
University of Bayreuth as a Humboldt Research Fellow. The
author thanks Alexander von Humboldt Foundation for the
financial support and the University of Bayreuth for the
hospitality. I am especially thankful to Prof. Th. Peternell.

The paper is an enlarged version of the talk, where some
details and an explanation of the technique of hypertangent
divisors, making the result of the paper [P] slightly
stronger, are added. The present paper was written during 
my stay at Max-Planck-Institut f\" ur Mathematik in Bonn.
I am very grateful to the Institute for the financial
support, stimulating atmosphere and hospitality.
\vspace{1cm}


{\large\bf 1. The Noether theorem}
\vspace{0.5cm}

\noindent
In [N] Max Noether published his famous theorem on the
(two-dimensional) Cremona group: the group of birational
self-maps of the (complex) projective plane
$$
\mathop{\rm Bir} {\mathbb P}^2=\mathop{\rm Cr} {\mathbb P}^2=
\mathop{\rm Aut} {\mathbb C}(s,t)=\{\chi\colon {\mathbb P}^2 -\,-\,\to {\mathbb P}^2\}
$$
is generated by the group of projective automorphisms
$\mathop{\rm Aut} {\mathbb P}^2$ and a single 
quadratic Cremona transformation which in a
suitable coordinate system takes the form
\begin{equation}
\label{1}
\tau\colon (x_0:x_1:x_2)\mapsto (x_1x_2:x_0x_2:x_0x_1).
\end{equation}
His argument went as follows. Take an arbitrary birational 
self-map $\chi$ of ${\mathbb P}^2$ and consider the strict
transform of the linear system of lines via $\chi^{-1}$:
$$
\begin{array}{ccc}
{\mathbb P}^2    &   \stackrel{\chi}{-\,-\,\to}  &   {\mathbb P}^2  
\\  \\
\left\{\textstyle
\begin{array}{c}
\mbox{curves} \\
\mbox{of degree} \\
n\geq 1
\end{array}\right\}
  &    \longleftarrow  & \{\mbox{lines}\}  
\\  \\
\|  &   &    
\\  \\
\left\{\textstyle
\begin{array}{c}
\mbox{the linear} \\
\mbox{system}\,\, |\chi|
\end{array}\right\} 
\end{array}
$$
The moving (that is, free from fixed components) linear
system $|\chi|$ becomes naturally the main subject of study.
One has the following obvious alternative:

\begin{itemize}

\item  either $n=1$, in which case
$\chi\in\mathop{\rm Aut}{\mathbb P}^2$ is regular,

\item  or $n\geq 2$, in which case $\chi$ is a birational
map in the proper sense; in particular, the linear
system $|\chi|$ has base points.

\end{itemize}

Assume that the second case holds. Since the curves in the
linear system $|\chi|$ are rational and the {\it free intersection} (that is, the intersection outside the base
locus) is equal to one, one can deduce that there exist at
least three distinct points $o_1,o_2,o_3$ of the linear
system $|\chi|$ satisfying the {\it Noether inequality}:
$$
\sum^3_{i=1}\mathop{\rm mult}\nolimits_{o_i}|\chi|>n.
$$
If all three points $o_i$ lie on ${\mathbb P}^2$ (that is, there are
no infinitely near points among them), then take the
standard Cremona transformation $\tau$ (\ref{1}) where
$o_1,o_2,o_3$ are assumed to be the points
$$
(1,0,0),\, (0,1,0,)\,\, \mbox{and}\,\, (0,0,1),
$$
respectively. It is easy to compute that the linear system
$|\chi\circ\tau|$ (which is the strict transform of $|\chi|$
via $\tau$, or the strict transform of the linear system of lines on ${\mathbb P}^2$ via the composition $\chi\circ\tau$) is a
moving linear system of plane curves of degree
$$
2n-\sum^3_{i=1}\mathop{\rm mult}\nolimits_{o_i}|\chi|<n.
$$
In other words, taking the composition with a quadratic
transformation (which are all conjugate with each other
by a projective automorphism), one can decrease the
degree $n\geq 1$. Thus (assuming that at each step there
are no infinitely near points among $o_i$'s) we get a 
decomposition of $\chi$ into a product of quadratic
transformations:
$$
\chi=\tau_1\circ\dots\circ\tau_N\circ\alpha=
\alpha_1\circ\tau\circ\alpha_2\circ\tau\circ\dots
\circ\tau\circ\alpha_{N+1},
$$
where $\tau$ is the standard involution (\ref{1}) and
$\alpha,\alpha_i$ are all projective automorphisms.

There is an immense literature on the Noether theorem
and Cremona transformations, see, for instance [H]
(the book is to be soon re-published with an
explanatory introduction written by V.A.Iskovskikh and 
M.Reid). Here we are interested only in the principal 
ingredients of Noether's argument. These are:

\begin{itemize}

\item the invariant $n\geq 1$ (the degree of curves
in the linear system $|\chi|$),

\item  ``maximal'' triples of base points (that is,
the triples satisfying the Noether inequality (\ref{1})),

\item the ``untwisting'' procedure (decreasing $n$ and
thus ``simplifying'' $\chi$).

\end{itemize}

{\bf Remark.} One can well imagine that the untwisting
procedure may not be uniquely determined. This is the case,
when there are more than one ``maximal'' triples, so that
we can decrease the degree $n$ taking the composition
with various quadratic involutions. This naturally
leads to {\it relations} between the generators of
$\mathop{\rm Bir} {\mathbb P}^2$. They were first described by M.Gizatullin in [G]. Later the argument was radically simplified by V.A.Iskovskikh [I4].
\vspace{1cm}


{\large\bf 2. Fano's work}
\vspace{0.5cm}

\noindent
At the beginning of the XXth century Fano started his work
in three-dimensional birational geometry. He was an absolute
pioneer in the field. His work lasted for about 50 years
and, apart from its great mathematical value, presents an
example of an equally great courage and inner strength.

Fano started with an attempt to reproduce Noether's
argument in dimension three. His first object of study was
the famous three-dimensional quartic 
$V=V_4\subset{\mathbb P}^4$. His
investigation went as follows. Take a birational self-map
$\chi\in\mathop{\rm Bir} V$ and look at the strict transform of the
linear system of hyperplane sections of $V$ via $\chi^{-1}$:
$$
\begin{array}{rcl}
V    &   \stackrel{\chi}{-\,-\,\to}  &   V
\\
\left\{\textstyle
\begin{array}{c}
\mbox{a linear system} \\
\mbox{of divisors}\,\, |\chi|  \\
\mbox{cut out on}\,\, V\,\, \mbox{by}  \\
\mbox{hypersurfaces} \\
\mbox{of degree}\,\, n\geq 1 \\
\end{array}\right\}
  &    \longleftarrow  & 
\left\{\textstyle
\begin{array}{c}
\mbox{a linear system} \\
\mbox{of hyperplane} \\
\mbox{sections}
\end{array}\right\} 
\end{array}
$$
Now, similar to the two-dimensional case, we get an
obvious alternative:

\begin{itemize}

\item  either $n=1$, in which case $\chi\in\mathop{\rm Aut} V$
is regular,

\item or $n\geq 2$, in which case $\chi$ is a birational
map in the proper sense; in particular, the linear
system $|\chi|$ has a non-empty base locus.

\end{itemize}

According to the scheme of Noether's arguments, the
next step to be made is finding a subscheme of high
multiplicity in the base scheme of the linear system
$|\chi|$. And indeed, Fano asserted that if $n\geq 2$,
then one of the following possibilities holds:

\begin{itemize}

\item  there exists a curve $B\subset V$ such that
\begin{equation}
\label{2}
\mathop{\rm mult}\nolimits_B |\chi|>n,
\end{equation}

\item  there exists a point $x\in V$ such that
\begin{equation}
\label{3}
\mathop{\rm mult}\nolimits_x |\chi|>2n,
\end{equation}

\item something similar happens, reminding of the
unpleasant infinitely near case of the Noether theorem.
Here Fano does not give any formal description, just 
presents an example of what can take place: there is
a point $x\in V$ and an infinitely near line 
$L\subset E\cong{\mathbb P}^2$, where
$$
\varphi\colon \widetilde V \to V
$$
is the blowing up of $x$ with the exceptional divisor
$E$, such that
$$
\mathop{\rm mult}\nolimits_x|\chi|+\mathop{\rm mult}\nolimits_L \widetilde{|\chi|}> 3n,
$$
where $\widetilde{|\chi|}$ is the strict transform
of the linear system $|\chi|$ on $\widetilde V$.

\end{itemize}

Now Fano gives certain arguments, some of which are true
and some not, showing that these cases are impossible.
He concludes that the $n\geq 2$ case does not realize
and therefore $n=1$ is the only possible case. Thus
$$
\mathop{\rm Bir} V=\mathop{\rm Aut} V.
$$

Later Fano studied several other types of three-folds. One of his
most impressive claims concerns the complete intersection
$$ 
\begin{array}{rcl} \displaystyle
V= & V_{2\cdot 3} &
\subset {\mathbb P}^{5}\\ \displaystyle
& \parallel & \\ \displaystyle &
F_2\cap F_3 &
\end{array}
$$
of a quadric and a cubic in ${\mathbb P}^5$. Starting as above,
$$
\begin{array}{rcl}
V    &   \stackrel{\chi}{-\,-\,\to}  &   V
\\
|nH|\supset |\chi|
  &    \longleftarrow  & 
\left\{\textstyle
\begin{array}{c}
\mbox{the linear system} \\
\mbox{of hyperplane} \\
\mbox{sections}
\end{array}\right\} 
\end{array}
$$
($H\in\mathop{\rm Pic} V$ is the class of a hyperplane section,
$\mathop{\rm Pic} V={\mathbb Z} H$), Fano discovers that for certain special
subvarieties his inequality (\ref{2}) can be satisfied.
This is the case, when $B=L\subset V$ is a line. It is 
easy to see that the projection $V-\,-\,\to {\mathbb P}^3$
from the line $L$ is a dominant rational map of degree 2.
Thus there exists a Galois involution
$$
\tau_L\in \mathop{\rm Bir} V,
$$
permuting points in a general fiber. One computes easily
that
$$
|\tau_L|\subset |4H|\quad\mbox{and}\quad
\mathop{\rm mult}\nolimits_L|\tau_L|=5,
$$
so that the linear system $|\tau_L|$ realizes the
inequality (\ref{2}). Now if
$$
\mathop{\rm mult}\nolimits_L|\chi|>n, \quad
|\chi|\subset |nH|,
$$
then one can check that the linear system
$|\chi\circ\tau_L|$ is cut out on $V$
by hypersurfaces of degree
$$
4n-3\mathop{\rm mult}\nolimits_L|\chi|<n,
$$
which gives the necessary ``untwisting'' procedure for
$\chi$. The analogy to Noether's arguments is now
complete. Let us once again look at the general scheme
of Fano's arguments. They consist of the following components:

\begin{itemize}

\item the invariant $n$ ($|\chi|\subset |nH|$),

\item  existence of ``maximal'' curves or points
(or something similar), satisfying the Fano inequalities
(\ref{2},\ref{3}),

\item  either excluding or untwisting the maximal curve
or point found at the previous step (decreasing $n$ and
``simplifying'' $\chi$).

\end{itemize}

It should be added that the untwisting procedure is not 
always uniquely determined, that is, Fano inequalities are
sometimes satisfied for a few various subvarieties, e.g.
two different lines on $V_{2\cdot 3}$ (when the
corresponding plane is contained in the quadric
$F_2$). This 
naturally leads to relations between generators of
$\mathop{\rm Bir} V$. For $V_{2\cdot 3}$ they were described by
V.A.Iskovskikh around 1975, see [I3] and a detailed
exposition in [IP].
\vspace{1cm}


{\large\bf 3. The theorem of V.A.Iskovskikh and Yu.I.Manin}
\vspace{0.5cm}

\noindent
The modern birational geometry of three-dimensional
varieties started in 1970 with two major breakthroughs:
the theorem of H.Clemens and Ph.Griffiths on the
three-dimensional cubic [CG] and the theorem of V.A.Iskovskikh and Yu.I.Manin on the three-dimensional
quartic [IM]. In the latter paper Fano's ideas were
developed into a rigorous and powerful theory, which made
it possible to begin a systematic study of explicit 
birational geometry of three-folds. This success was to a
considerable degree prepared by the earlier papers of Yu.I.Manin on surfaces over non-closed fields: in [M1, M2]
all the principal technical components of the method
of maximal singularities were already present, including
the crucial construction of the graph, associated
with a finite sequence of blow ups.

Let us reproduce briefly the arguments of [IM]. Fix a
smooth quartic $V=V_4\subset{\mathbb P}^4$ and consider, as usual,
the strict transform of the linear system of hyperplane sections with respect to $\chi^{-1}$, where
$\chi\colon V -\,-\,\to V'_4$ is a birational map onto another smooth quartic:
$$
\begin{array}{rcl}
V    &   \stackrel{\chi}{-\,-\,\to}  &   V'
\\
|nH|\supset |\chi|
  &    \longleftarrow  & 
\left\{\textstyle
\begin{array}{c}
\mbox{the linear system} \\
\mbox{of hyperplane} \\
\mbox{sections}
\end{array}\right\} 
\end{array}
$$
As usual, we come to the familiar alternative:

\begin{itemize}

\item  either $n=1$, in which case $\chi\colon V\to V'$
is an isomorphism (hence a projective isomorphism),

\item or $n\geq 2$, in which case $\chi$ is a birational
map in the proper sense; in particular, the base
subscheme of the linear system $|\chi|$ is non-empty.

\end{itemize}

Assume that $n\geq 2$.
\vspace{0.3cm}

{\bf Proposition 3.1.} {\it There exists a geometric
discrete valuation $\nu$ on $V$ (here ``geometric'' means
``realizable by a prime divisor $E$ on some model 
$\widetilde V$ of the
field ${\mathbb C}(V)$'') such that the 
inequality
\begin{equation}
\label{4}
\nu(\Sigma)>n\cdot\mathop{\rm discrepancy}(\nu)
\end{equation}
holds, where $\nu(|\chi|)=\nu(D)$
for a general divisor $D\in|\chi|$.}
\vspace{0.3cm}

 The discrete
valuations $\nu$, satisfying (\ref{4}), are called
{\it maximal singularities} of the linear system $|\chi|$.
The inequality (\ref{4}) is called the
{\it Noether-Fano inequality}.

Now [IM] shows that a moving linear system $|\chi|$ on
$V$ cannot have a maximal singularity. The hardest case
is when the centre of the maximal singularity $\nu$ is a point: $\mathop{\rm centre}(E)=x\in V$. In this case take 
two general divisors $D_1,D_2\in|\chi|$ and consider the
cycle of scheme-theoretic intersection
$$
Z=(D_1\circ D_2).
$$
It is an effective curve on $V$. The crucial fact is
given by
\vspace{0.3cm}

{\bf Proposition 3.2.} {\it The following inequality holds}
$$
\mathop{\rm mult}\nolimits_x Z>4n^2.
$$
\vspace{0.3cm}

Since $\mathop{\rm deg} Z=4n^2$, this gives a contradiction. Thus we obtain
\vspace{0.5cm}

{\bf Theorem 3.1.[IM]} {\it Any birational map between smooth
three-dimensional quartics is a projective isomorphism.}
\vspace{0.5cm}

However this very argument gives immediately a much
stronger claim! For instance, let us describe birational
maps
$$
\begin{array}{cccl}
V  & \stackrel{\chi}{-\,-\,\to}  & V'  &   \\
   &                        &  \downarrow   &  \pi'  \\
   &                        &  S'  &   
\end{array}
$$
from $V$ to conic bundles $V'/S'$ (or, in a slightly
different language, describe the {\it conic bundle
structures} on $V$). Let us construct a moving linear
system $\Sigma'$ on $V'$ in the following way:
$$
\begin{array}{rccl}
   &  V'   &   \Sigma'   &   \\  \\
\pi'  & \downarrow  &  \uparrow & \mbox{pull back via}
\,\,\pi'  \\
   &  S'  &  \left\{
\begin{array}{c}
\mbox{moving linear}\\
\mbox{system of curves} \\
\mbox{on the surface}\,\, S'
\end{array}
\right\}
\end{array}
$$
Consider the strict transform $\Sigma\subset |nH|$
of $\Sigma'$ on $V$ via $\chi$. Then it is easy to prove
that a maximal singularity exists {\it always}, now
irrespective of the value of $n\geq 1$. However, we
know that existence of a maximal singularity leads to a
contradiction. Therefore, the birational map $\chi$
simply cannot exist!

What was actually proved in [IM],
can be formulated as follows:
\vspace{0.5cm}

{\it a smooth three-dimensional quartic $V\subset{\mathbb P}^4$

\begin{itemize}

\item cannot be fibered into rational curves by a
rational map,

\item cannot be fibered into rational surfaces by a
rational map,

\item if $\chi\colon V-\,-\,\to V'$ is a birational map,
where $V'$ is a ${\mathbb Q}$-Fano threefold, then
$\chi\colon V\to V'$ is an isomorphism.

\end{itemize}
}

Speaking the modern language, we express all this by
saying that the quartic is {\it birationally superrigid}.
\vspace{1cm}


{\large\bf 4. The method of maximal singularities}
\vspace{0.5cm}

\noindent
In all the procedures described above a certain integral
parameter was involved --- namely, the ``degree'' $n$  of
the linear system $\Sigma$, defining the birational map
$\chi$ under consideration. All the above examples dealt
with Fano varieties $V$ such that $\mathop{\rm Pic} V\cong {\mathbb Z}$, so
that the integer $n$ meant just the class of $\Sigma$
in $\mathop{\rm Pic} V$. However, this extremely important number
has a more general invariant meaning, which we describe
now.

Let $X$ be a uniruled ${\mathbb Q}$-Gorenstein variety
with terminal singularities. This assumption implies
that the canonical class $K_X$ is negative on some
family of (generically irreducible) curves sweeping out $X$.
Therefore for any divisor $D$ the following number is
finite:
$$
c(D,X)=\sup\{b/a|b,a\in{\mathbb Z}_+\setminus \{0\},
|aD+bK_X|\neq\emptyset\}.
$$
It is called the {\it threshold of canonical adjunction}
of the divisor $D$. Sometimes we omit $X$ and write
simply $c(D)$ or $c(\Sigma)$ for $D\in \Sigma$ moving
in a linear system.

Now we can describe the general scheme of the method of
maximal singularities. Let us fix a uniruled variety
$V$ with ${\mathbb Q}$-factorial terminal singularities
and another variety $V'$ in this class. Let us assume that
$V'$ is birational to $V$. The aim of the method is to
give a complete description of all possible birational
correspondences between $V$ and $V'$.

We start as usual with the following diagram:

$$
\begin{array}{rcccl}
   &  V  &  \stackrel{\chi}{-\,-\,\to}  &
 V'  &    \\  \\
\begin{array}{c}
\mbox{moving linear} \\
\mbox{system}
\end{array}
&  \Sigma  & \longleftarrow & \Sigma' & 
\begin{array}{c}
\mbox{moving linear} \\
\mbox{system}
\end{array}  \\  \\
   &  c(\Sigma) & ? & c(\Sigma') &

\end{array}
$$
The linear system $\Sigma'$ is fixed throughout the whole
argument. The thresholds are taken with respect to the
varieties $V$, $V'$. The ? sign means that we do not know,
which inequality is true: ``$\leq$'' or ``$>$''.

Now we get the alternative:
\begin{itemize}

\item either $?$ is $\leq$, in which case we stop. It is
presumed that when we have the inequality
\begin{equation}
\label{5}
c(\Sigma)\leq c(\Sigma'),
\end{equation}
then ``we can say everything'' about the map $\chi$.
In real life, sometimes this is the case, sometimes not.
But in any case this inequality completely reduces the
birational problem to a biregular one, since the family
of linear systems $\Sigma$ is bounded (in many cases
(\ref{5}) implies that it is actually empty or $\Sigma$
is unique, as above).

\item  or $?$ is $>$, in which case we proceed further
as follows.

\end{itemize}
\vspace{0.3cm}

{\bf Proposition 4.1.} {\it There exists a geometric
discrete valuation $\nu=\nu_E$ on $V$ such that
the inequality 
\begin{equation}
\label{6}
\nu(\Sigma)>c(\Sigma)\cdot a(E)
\end{equation}
holds.
}
\vspace{0.3cm}

The inequality (\ref{6}) is called {\it the Noether-Fano inequality}. The discrete valuation $\nu$ is called a
{\it maximal singularity} of the linear system $\Sigma$.

The word ``geometric'' means, as we have mentioned above, 
that there exists a birational
morphism $\varphi\colon\widetilde V \to V$ with $\widetilde V$ smooth
such that $\nu=\nu_E$ for some prime divisor 
$E\subset\widetilde V$. In (\ref{6}) $\nu_E(\Sigma)$ means the
multiplicity of a general divisor $D\in\Sigma$ at $E$
and $a(E)$ means the discrepancy of $E$.

Now for {\it each} geometric discrete valuation $E$ the
following work should be performed:

\begin{itemize}

\item either $E$ can be excluded as a possible maximal
singularity: there is no {\it moving} linear system
$\Sigma$ satisfying the Noether-Fano inequality
(\ref{6}) for $E$,

\item  or $E$ should be untwisted. The untwisting means
that we find a birational self-map
$$
\chi^*_E\in \mathop{\rm Bir} V
$$
such that we get the following diagram:

$$
\begin{array}{ccc}
V  &  \stackrel{\chi^*_E}{-\,-\,\to}  &  V  \\  \\
\Sigma^*  &  \longleftarrow & \Sigma    \\  \\
c(\Sigma^*)  &  <  &  c(\Sigma).
\end{array}
$$
Replacing $\chi$ by
$$
\chi\circ\chi^*_E\colon V -\,-\,\to V',
$$
we go back to the beginning of the procedure
(compare $c(\Sigma^*)$ and $c(\Sigma')$ and so on).

\end{itemize}
\vspace{1cm}


{\large\bf 5. Birationally rigid varieties}
\vspace{0.5cm}

\noindent
Roughly speaking $V$ is said to be
{\it birationally rigid}, if the above-described procedure
works for $V$, that is, in a finite number of steps
we obtain the desired inequality (\ref{5}).
\vspace{0.1cm}

{\bf Definition 5.1.} $V$ is said to be birationally rigid,
if for {\it any} $V'$, {\it any} birational map
$\chi\colon V-\,-\,\to V'$ and {\it any} moving linear system $\Sigma'$
on $V'$ there exists a birational self-map
$\chi^*\in \mathop{\rm Bir} V$ such that the following diagram holds:

$$
\begin{array}{ccccc}
V  & \stackrel{\chi^*}{-\,-\,\to}
&  V  &  \stackrel{\chi}{-\,-\,\to} & V'\\  \\
\Sigma^*  & \longleftarrow & \Sigma  &
\longleftarrow & \Sigma'  \\  \\
c(\Sigma^*)  &  &  \leq &  &  c(\Sigma')
\end{array}
$$
\vspace{0.1cm}

\noindent
The birational self-map $\chi^*$ is a composition of
elementary untwisting maps described above:
$$
\chi^*=\chi^*_{E_1}\circ\dots\circ \chi^*_{E_N}.
$$
However, it turns out that for many (hopefully, for
``majority'' of) Fano varieties the untwisting procedure
is trivial.
\vspace{0.1cm}

{\bf Definition 5.2.} $V$ is said to be birationally superrigid,
if for {\it any} $V'$, {\it any} birational map
$V-\,-\,\to V'$ and {\it any} moving linear system $\Sigma'$
on $V'$ the following diagram holds:

$$
\begin{array}{ccc}
V  &  \stackrel{\chi}{-\,-\,\to}  &  V'  \\  \\
\Sigma  &  \longleftarrow & \Sigma'    \\  \\
c(\Sigma)  &  <  &  c(\Sigma').
\end{array}
$$
\vspace{0.1cm}

In other words, a birationally rigid variety is
superrigid, if we may always take 
$\chi^*=\mathop{\rm id}_V$. In the sense of the given definitions, the smooth quartic $V_4\subset{\mathbb P}^4$ is
superrigid and the smooth complete intersection
$V_{2\cdot 3}\subset {\mathbb P}^5$ is rigid (for the latter
case, the proof has been so far produced for a generic
member of the family only, see [I3,P2,IP]).

Immediate geometric implications of birational 
(super)rigidity, which actually determined the very choice
of this word combination, are collected below.
\vspace{0.3cm}

{\bf Proposition 5.1.} {\it Let $V$ be a smooth Fano variety
with $\mathop{\rm Pic} V={\mathbb Z} K_V$. If $V$ is birationally rigid, then

\begin{itemize}

\item  $V$ cannot be fibered into rationally connected
(or uniruled) varieties by a non-trivial rational map;
that is, the following diagram is impossible
$$
\begin{array}{cccl}
V  &  {-\,-\,\to} & V'  &  \\
  &  &  \downarrow & \mbox{\rm uniruled fibers} \\
  &  &  S'  &   
\end{array}
$$
with $\mathop{\rm dim}S'\geq 1$,

\item  if $\chi\colon V-\,-\,\to V'$ is a birational map onto
a ${\mathbb Q}$-Fano variety with 
$\mathop{\rm rk} \mathop{\rm Pic} V'=1$, then
$$
V\cong V'
$$
(although $\chi$ itself may be {\it not} a biregular map).
If, moreover, $V$ is superrigid, then $\chi$ itself
is an isomorphism. In particular, in the superrigid case
the groups of birational
and biregular self-maps coincide,
\begin{equation}
\label{7}
\mathop{\rm Bir} V=\mathop{\rm Aut} V.
\end{equation}
Conversely, if $V$ is rigid and (\ref{7}) holds, then by
definition of rigidity it is clear that $V$ is superrigid.

\end{itemize}
}
\vspace{0.3cm}

The known examples motivate the following
\vspace{0.5cm}

{\bf Conjecture 5.1.} {\it Let $V$ be a smooth Fano variety
of dimension $\mathop{\rm dim} V\geq 4$ with
$\mathop{\rm Pic} V={\mathbb Z} K_V$. Then $V$ is birationally rigid. If
$\mathop{\rm dim} V\geq 5$, then $V$ is superrigid.
}
\vspace{0.5cm}

The biggest class of higher-dimensional Fano varieties,
supporting this conjecture, is given by the following
\vspace{0.5cm}

{\bf Theorem 5.1. [P5,P8]} {\it Let
$$ 
\begin{array}{rcl} \displaystyle
V= & V_{d_1\cdot d_2\cdot \dots \cdot d_k} &
\subset {\mathbb P}^{M+k}\\ \displaystyle
& \parallel & \\ \displaystyle &
F_1\cap F_2\cap \dots \cap F_k &
\end{array}
$$
be a sufficiently general (in the sense of Zariski
topology) Fano complete intersection of the type
$d_1\cdot d_2\cdot \dots \cdot d_k$, where 
$d_1+\dots +d_k=M+k$, $M\geq 4$ and 
$2k<M=\mathop{\rm dim}V$. Then $V$ is birationally
superrigid.
}
\vspace{0.5cm}

More examples of superrigid Fano varieties are given by
smooth complete intersections in weighted projective
spaces, see [P6].

{\bf Remark.} As A.Beauville kindly informed the author,
in dimension 4 there is still a gap between rigidity and
superrigidity. His example is a $2\cdot 2\cdot 3$ smooth
complete intersection in ${\mathbb P}^7$, containing a 
two-dimensional plane $P\cong{\mathbb P}^2$. There are no reasons 
to doubt their rigidity (although it has not yet been
proved). However, for these special complete intersections
$$
\mathop{\rm Bir} V\neq \mathop{\rm Aut} V,
$$
since the projection from the plane $P$,
$$
\pi\colon V -\,-\,\to {\mathbb P}^4,
$$
is of degree two and therefore generates a Galois
involution $\tau\in \mathop{\rm Bir} V$ which is not
regular. Therefore such varieties are not
superrigid. In dimension $\geq 5$ constructions of this
type are not possible.

Most of the proofs of birational superrigidity
make use of the following sufficient condition.
\vspace{0.5cm}

{\bf Theorem 5.2.} {\it Let $X$ be a smooth Fano variety with
$\mathop{\rm Pic} X={\mathbb Z} K_X$. Assume that for any irreducible
subvariety $Y\subset X$ of codimension two
the following two properties are satisfied:} \par
\vspace{0.1cm}

(i) $\mathop{\rm mult}\nolimits_Y\Sigma\leq n$
{\it for any linear system $\Sigma\subset |-nK_X|$
without fixed components};\par
\vspace{0.1cm}

(ii) {\it the inequality}
\begin{equation}
\label{c1}
\mathop{\rm mult}\nolimits_x Y\leq
\frac{\displaystyle 4}{\displaystyle
\mathop{\rm deg} X}
\mathop{\rm deg} Y
\end{equation} 
{\it  holds for any point} $x\in Y$, {\it where 
$$
\mathop{\rm deg} X=(-K_X)^{\mathop{\rm dim} X},\quad
\mathop{\rm deg} Y=(Y\cdot (-K_X)^{\mathop{\rm dim} Y})
$$
and $\mathop{\rm mult}\nolimits_Y\Sigma$ means multiplicity
of a general divisor $D\in\Sigma$ along $Y$. 
\vspace{0.1cm}

Then the variety $X$ is birationally superrigid.} 
\vspace{0.5cm}

The strongest technique which makes it possible to check
the (principal) condition (ii) of this criterion is
that of {\it hypertangent divisors}. Although at the moment
an alternative method, based on the connectedness principle
of Shokurov and Koll\' ar [Sh,K] (suggested by Corti [C2]
and later used in [CM] and [P10]), gains momentum, the
older argument by hypertangent divisors is still
working  better. For the reader to get the idea of this
technique, we give here a proof for Fano hypersurfaces
$V=V_M\subset{\mathbb P}^M$. Birational superrigidity of any smooth
hypersurface has already been proved in [P10]. Nevertheless
we give here an argument which is based on the paper [P5],
slightly sharpening the result.
\vspace{0.3cm}

{\bf Proposition 5.2.} {\it Let $x\in V=V_M\subset{\mathbb P}^M$ be
a point such that there are but finitely many lines
$L\subset V$ passing through $x$. Then for any irreducible
subvariety $Y\subset V$ of codimension two the estimate
\begin{equation}
\label{8}
\frac{\mathop{\rm mult}\nolimits_x}{\mathop{\rm deg}} Y\leq \frac{\displaystyle 4}{\displaystyle \mathop{\rm deg} V}=
\frac{4}{M}
\end{equation}
holds.
}
\vspace{0.3cm}

{\bf Proof.} Let $(z_1,\dots,z_M)$ be a system of affine
coordinates on ${\mathbb P}={\mathbb P}^M$ with the origin at $x$. Write
down the equation of the hypersurface $V$:
$$
f=q_1+q_2+\dots+q_M,
$$
where $q_i$ are homogeneous of degree $i$ in $z_*$.
Note that the lines through $x$ on the hypersurface $V$
are given by the system of equations
\begin{equation}
\label{9}
q_1=q_2=\dots=q_M=0.
\end{equation}
Therefore the set (\ref{9}) of common zeros is of
dimension at most one. Denote by
$$
f_i=q_1+\dots+q_i
$$
the truncated polynomials. It is clear that in the
affine open set ${\mathbb A}={\mathbb A}^M_{(z_1,\dots,z_M)}\subset{\mathbb P}$
the algebraic set
$$
f_1=f_2=\dots=f_M=0
$$
is the same as (\ref{9}), therefore it is of dimension
at most one. This implies, in its turn, that the
algebraic set
$$
f_1|_{{\mathbb A}\cap V}=f_2|_{{\mathbb A}\cap V}=\dots
=f_{M-1}|_{{\mathbb A}\cap V}=0
$$
on the affine part of the hypersurface $V$ is also of
dimension at most one: scheme-theoretically it is the
same as (\ref{9}), supported on the union of lines on $V$
through $x$.

Let us look at the divisors
$$
D_i=\overline{
\{f_i|_{{\mathbb A}\cap V}=0\}
},
$$
$i=1,\dots,M-1$. We call them {\it hypertangent divisors}:
if $H\in\mathop{\rm Pic} V$ is the class of a hyperplane section, then
clearly
$$
D_i\in |iH|
$$
and
$$
\mathop{\rm mult}\nolimits_x D_i\geq i+1,
$$
since in the affine part of $V$
$$
D_i|_{{\mathbb A}\cap V}=\{
(q_{i+1}+\dots+q_M)|_V=0\}.
$$
Now by assumption
$$
\mathop{\rm dim}\nolimits_x(D_1\cap\dots\cap D_{M-1})\leq 1,
$$
where $\mathop{\rm dim}\nolimits_x$ means the dimension in
a neighborhood of the point $x$. It is easy to see that for
a given subvariety $Y\ni x$ of codimension two there is
a set of $(M-4)$ divisors
$$
\{D_i\,|\, i\in{\cal I}\}\subset
\{D_1,\dots,D_{M-1}\}
$$
such that
$$
\mathop{\rm dim}\nolimits_x 
\left(
Y\cap\mathop{\bigcap}\limits_{i\in {\cal I}}D_i
\right)=1.
$$
Now let us order the set ${\cal I}$ somehow, so that
$$
\{D_i\,|\, i\in{\cal I}\}=\{R_1,R_2,\dots,R_{M-4}\}.
$$
It is easy to construct by induction on
$i\in\{0,\dots,M-4\}$ the sequence of irreducible 
subvarieties
$$
Y_0=Y,Y_1,\dots,Y_{M-4},
$$
such that

\begin{itemize}

\item $Y_{i+1}\subset Y_i$, 
$\mathop{\rm codim} Y_i=i+2$ (the codimension is
taken with respect to $V$);

\item  $Y_i\not\subset R_{i+1}$, so that
$(Y_i\circ R_{i+1})$ is an effective algebraic cycle
on $V$, $Y_{i+1}$ is one of its irreducible components;

\item the following estimate holds:
$$
\frac{\mathop{\rm mult}\nolimits_x}{\mathop{\rm deg}} Y_{i+1}\geq \frac{\mathop{\rm mult}\nolimits_x}{\mathop{\rm deg}} Y_i\cdot
\frac{\displaystyle \mathop{\rm mult}\nolimits_x R_{i+1}}{a_{i+1}},
$$
where $R_j\in |a_j H|$, $j=1,\dots,M-4$.

\end{itemize}

Now $C=Y_{M-4}$ is an irreducible curve on $V$,
satisfying the inequality
\begin{equation}
\label{10}
\frac{\mathop{\rm mult}\nolimits_x}{\mathop{\rm deg}} C\geq \frac{\mathop{\rm mult}\nolimits_x}{\mathop{\rm deg}} Y
\cdot \frac54 \cdot \cdots \cdot \frac{M}{M-1}.
\end{equation}
(If the set ${\cal I}\neq \{4,\dots,M-1\}$, then
the estimate is better: we take the worst possible
case.)
Making the obvious cancellations and taking into
consideration that the left-hand side of
(\ref{10}) cannot exceed 1, we obtain the desired
estimate (\ref{8}).
\vspace{0.3cm}

{\bf Remark.} When we define informally birationally
rigid varieties as those for which the method of maximal
singularities works, one may naturally ask, what happens
if it does not. There is an answer in dimension three.
It is given by the Sarkisov program, developed by
V.G.Sarkisov (see [S3]) and completely proved by Corti
in [C1]. The answer is, that when it is possible
neither to exclude nor to untwist a maximal singularity, it should be eliminated by a link to another Mori fiber space.
We do not touch these points in the present paper.
See [S3,R,C1,C2,CR,CM] for the details.
\vspace{1cm}


{\large\bf 6. Singular Fano varieties}
\vspace{0.5cm}

\noindent
Up to this moment, all our examples dealt with smooth
Fano varieties. Here we give the most interesting cases
of birationally rigid Fano varieties with isolated terminal
singularities. The oldest example of a birationally 
rigid singular Fano 3-fold is the three-dimensional
quartic with a unique non-degenerate double point [P1].

Let $x\in V=V_4\subset {\mathbb P}^4$ be the singularity. There
are 24 lines on $V$ passing through $x$; denote them by
$L_1,\dots L_{24}$. With the point $x\in V$ a birational
involution $\tau\in \mathop{\rm Bir} V$ is naturally associated: the
projection from $x$
$$
\begin{array}{rcccl}
   &   V   &   \subset  &  {\mathbb P}^4  &   \\
\mbox{rational map} & |  &  &  |  & \mbox{rational map} \\
\mbox{of degree 2}  & |  &  &  |  &  \mbox{with the fiber}\,\, {\mathbb P}^1  \\
          & \downarrow  &  &  \downarrow  &  \\
    &  {\mathbb P}^3  &  =  &  {\mathbb P}^3  & 
\end{array}
$$
determines the Galois involution $\tau$ of $V$ over
${\mathbb P}^3$.

Let $L=L_i$ be a line through $x$. Look at the projecttion
from the line $L$:

$$
\begin{array}{rcccl}
   &   V   &   \subset  &  {\mathbb P}^4  &   \\
\mbox{rational map} & |  &  &  |  & \mbox{rational map} \\
\mbox{with the fibers ---}  & |  &  &  |  & \mbox{with the 
fiber}\,\, {\mathbb P}^2
          \\
\mbox{cubic curves}  & \downarrow  &  &  \downarrow  &  \\
    &  {\mathbb P}^2  &  =  &  {\mathbb P}^2.  & 
\end{array}
$$
Thus $V$ is fibered over ${\mathbb P}^2$ into elliptic curves.
Since $V$ is singular at $x$, all the cubic curves pass through $x$, which means that the fibration $V/{\mathbb P}^2$
has a section. Taking it for the zero of a group law
on a generic fiber, we get the birational involution
$$
\tau_i=\tau_L\in \mathop{\rm Bir} V,
$$
the reflection from zero on a generic fiber. Set
$\tau_0=\tau$. Let $B(V)$ be the subgroup of
$\mathop{\rm Bir} V$, generated by the involutions
$\tau_0,\tau_1,\dots,\tau_{25}$. 
\vspace{0.5cm}

{\bf Theorem 6.1. [P1]} {\it  {\rm (i)} $V$ is birationally rigid.

\noindent  {\rm (ii)} The group $B(V)$ is the free
product of 25 cyclic subgroups
$\langle \tau_i\rangle={\mathbb Z}/2{\mathbb Z}$, where
$i=0,1,\dots,25${\rm :}
$$
B(V)=\mathop{*}\limits^{24}_{i=0}\langle\tau_i\rangle.
$$

\noindent
{\rm (iii)} The subgroup $B(V)\subset\mathop{\rm Bir} V$ is normal and
the following exact sequence takes place:
$$
1\to B(V)\to \mathop{\rm Bir} V\to \mathop{\rm Aut} V\to 1.
$$
}
\vspace{0.5cm}

First proved in [P1], this theorem was later discussed
by Corti in [C2], where due to an application of the
connectedness principle of Shokurov and Koll\' ar [Sh,K]
the proof was simplified. The further study of singular
quartics was performed in [CM].

Theorem 6.1 can be generalized in higher dimensions [P9]. Let
$V=V_M\subset{\mathbb P}^M$, $M\geq 5$, be a sufficiently
general (for the precise conditions see [P9])
hypersurface with isolated terminal singularities. For
a singular point $x\in V$ we obviously have
$$
\mu_x=\mathop{\rm mult}\nolimits_x V\leq M-2,
$$
and if $\mu_x=M-2$, then the conditions of general
position imply that there is only one point with this
multiplicity. Let us define the integer
$$
\mu=\mathop{\rm max}\limits_{
x\in \mathop{\rm Sing} V} \{\mu_x\}\leq M-2.
$$
\vspace{0.5cm}

{\bf Theorem 6.2.} {\it {\rm (i)} Sufficiently general
variety $V$ is birationally rigid.

\noindent {\rm (ii)} If $\mu\leq M-3$, then $V$ is
superrigid.

\noindent {\rm (iii)} If $\mu=M-2$ and $\mu_x=\mu$
for the unique point $x$, then
$$
\mathop{\rm Bir} V=\langle \tau\rangle\cong
{\mathbb Z}/2{\mathbb Z},
$$
where $\tau\in\mathop{\rm Bir} V$ is the Galois involution
determined by the rational map
$$
\begin{array}{rcccl}
   &   V   &   \subset  &  {\mathbb P}^M  &   \\
\mbox{\rm rational map} & |  &  &  |  & \mbox{\rm projection} \\
\mbox{\rm of degree 2}  & |  &  &  |  &  
\mbox{\rm from}\,\, $x$         \\
          & \downarrow  &  &  \downarrow  &  \\
    &  {\mathbb P}^{M-1}  &  =  &  {\mathbb P}^{M-1}.  & 
\end{array}
$$
}
\vspace{0.5cm}

Of course, singular Fano varieties are much more numerous
in types than the smooth ones. The smooth quartics
generalize to 95 families of ${\mathbb Q}$-Fano
hypersurfaces
$$
V=V_d\subset {\mathbb P}(1,a_1,\dots,a_4),
$$
$a_1+a_2+a_3+a_4=d$. They all have terminal quotient
singularities but in a sense are closer by their
properties to smooth Fano varieties.
\vspace{0.5cm}

{\bf Theorem 6.3. [CPR]} {\it A general member $V$ of each
of 95 families is birationally rigid. The group of
birational self-maps is generated by finitely many
involutions.
}
\vspace{0.5cm}

For each of 95 families these involutions were explicitly
described in [CPR]. This paper is based on the classical
method of maximal singularities combined with the
Sarkisov  program [S3,R,C1].
\vspace{1cm}


{\large \bf 7. The relative version}
\vspace{0.5cm}

\noindent
So far we have been considering the absolute case, that is,
the case of Fano varieties. However, the world of 
rationally connected varieties is much bigger. If we
assume the predictions of the minimal model program,
each rationally connected variety is birational to a
fibration into Fano varieties over a base that, generally speaking, is not necessarily a point.

Here we give a very brief outline of the relative
version of the rigidity theory --- that is, rigidity 
theory of non-trivial fibrations. Similar to the
absolute case of Fano varieties, the starting point here
was formed by ``two-dimensional Fano fibrations'' over a
non-closed field, or simply speaking, surfaces with a
pencil of rational curves over a non-closed field. In the
papers of V.A.Iskovskikh [I1,I2] (which continued the work
started in the
papers of Yu.I.Manin on del Pezzo surfaces over non-closed
fields [M1,M2], see also [M3]) it was proved that under certain conditions
there is only one pencil of rational curves. This theorem
was the first relative rigidity result. It was necessary to
generalize these claims and, in the first place, the
technique of the proof to higher dimensions.

This breakthrough was made by V.G.Sarkisov in [S1,S2].
Let us consider smooth conic bundles of dimension $\geq 3$:
$$
\begin{array}{rrcccl}
   & &  V   &   \hookrightarrow  &  {\mathbb P}({\cal E})  &   \\  \\
\begin{array}{c}
\mbox{fibration}\\
\mbox{into conics}
\end{array}  & \pi  
 & \downarrow &  &  \downarrow  & 
\begin{array}{c}
\mbox{locally trivial}\\
{\mathbb P}^2-\mbox{fibration}
\end{array}
\\  \\
   &       &  S &  =  &  S  & 
\end{array}
$$
Here $\mathop{\rm dim}V\geq 3$,
$\mathop{\rm dim}S=\mathop{\rm dim}V-1\geq 2$,
${\cal E}$ is a locally trivial sheaf of rank 3 on $S$.
The points $x\in S$ over which the conic
$\pi^{-1}(x)\subset{\mathbb P}^2$ degenerates comprise
the {\it discriminant divisor} $D\subset S$.
Assume that $V/S$ is minimal (or {\it standard},
see [I3,S1,S2]) in the following sense:
$$
\mathop{\rm Pic} V={\mathbb Z} K_V \oplus \pi^* \mathop{\rm Pic} S.
$$
That is, $V/S$ is a Mori fiber space, see [C2].
In particular, $V/S$ has no sections. The main
question to be considered is whether $V$ has other
structures of a conic bundle or not. In other
words, let $\pi'\colon V'\to S'$ be another conic
bundle. Is an arbitrary birational map
$\chi\colon V -\,-\,\to V'$ fiber-wise or not?
Here is the diagram:
$$
\begin{array}{rcccl}
   &   V  &  \stackrel{\chi}{-\,-\,\to}  & V' &  \\
\pi  &  \downarrow &   &  \downarrow & \pi'  \\
   &  S  &  \stackrel{?}{-\,-\,\to} &  S',  &   
\end{array}
$$
where the ? sign means the question above.
\vspace{0.5cm}

{\bf Theorem 7.1. [S1,S2]} {\it If $|4K_S+D|\neq\emptyset$,
then $\chi$ is always fiber-wise: there exists a birational
map $\alpha\colon S -\,-\,\to S'$ such that
$$
\alpha\circ \pi=\pi'\circ\chi.
$$
}
\vspace{0.5cm}

There was no concept of birational rigidity in 1980. Now
we just say that $V/S$ is birationally rigid. In [S1,S2]
it was said just that the conic bundle structure, given by
definition, is unique.

Let $V/S$ be a non-trivial fibration, 
$\mathop{\rm dim} S\geq 1$, with rationally connected
(or just uniruled) fibers. Define the group of
{\it proper} birational self-maps setting
$$
\mathop{\rm Bir} (V/S)=\mathop{\rm Bir} F_{\eta}\subset \mathop{\rm Bir} V,
$$
where $F_{\eta}$ is the generic fiber, that is, the
variety $V$ considered over the field ${\mathbb C}(S)$. In
other words, birational self-maps from $\mathop{\rm Bir} (V/S)$
preserve the fibers.
\vspace{0.1cm}

{\bf Definition 7.1.} The fibration $V/S$ is
{\it birationally rigid},
if for {\it any} variety $V'$, {\it any} birational map
$V-\,-\,\to V'$ and {\it any} moving linear system $\Sigma'$
on $V'$ there exists a proper birational self-map
$\chi^*\in \mathop{\rm Bir} (V/S)$ such that the following diagram holds:

$$
\begin{array}{ccccc}
V  & \stackrel{\chi^*}{-\,-\,\to}
&  V  &  \stackrel{\chi}{-\,-\,\to} & V'\\  \\
\Sigma^*  & \longleftarrow & \Sigma  &
\longleftarrow & \Sigma'  \\  \\
c(\Sigma^*)  &  &  \leq &  &  c(\Sigma').
\end{array}
$$
\vspace{0.1cm}

The definition of superrigidity is word for word the same
as in the absolute case.
\vspace{0.3cm}

{\bf Proposition 7.1.} {\it Assume that $X/S$ is a Fano
fibration with $X$, $S$ smooth such that
$$
\mathop{\rm Pic} X={\mathbb Z}K_X\oplus
\pi^* \mathop{\rm Pic} S
$$
and for any effective class $D=mK_X+\pi^*T$ the
class $NT$ is effective on $S$ for some
$N\geq 1$. Assume furthermore
that $X/S$ is birationally rigid. Then for any
rationally connected fibration $X'/S'$ and any
birational map
$$
\chi\colon X -\,-\,\to X'
$$
(provided that such maps exist) there is a rational
dominant map
$$
\alpha\colon S -\,-\,\to S'
$$
making the following diagram commutative:
$$
\begin{array}{rcccl}
    &  X  &  \stackrel{\chi}{-\,-\,\to} & X' &  \\
\pi  &  \downarrow &    &  \downarrow  & \pi'  \\
    &  S   & \stackrel{\alpha}{-\,-\,\to} &  S'.  &
\end{array}
$$
}
\vspace{0.3cm}

{\bf Conjecture 7.1.} {\it If the fibration $V/S$ as above
is sufficiently twisted over the base $S$, then $V/S$ is
birationally rigid.}
\vspace{0.3cm}

We prefer not to be formal here about this conjecture. Instead of explaining what precisely is meant by the
twistedness assumption, we just give an example that
illustrates the point. This example has already been
generalized in higher dimensions [P7].

Let us consider three-folds fibered into cubic surfaces:
$$
\begin{array}{rrcccl}
   & &  V   &   \hookrightarrow  &  {\mathbb P}({\cal E})  &   \\  \\
\begin{array}{c}
\mbox{fibers are}\\
\mbox{cubic surfaces}
\end{array}  & \pi  
 & \downarrow &  &  \downarrow  & 
\begin{array}{c}
\mbox{locally trivial}\\
{\mathbb P}^3-\mbox{fibration}
\end{array}
\\   \\
   &       &  {\mathbb P}^1 &  =  &  {\mathbb P}^1.  & 
\end{array}
$$
Here $\mathop{\rm rk}{\cal E}=4$, 
$\mathop{\rm Pic} {\mathbb P}({\cal E})={\mathbb Z} L\oplus {\mathbb Z} G$,
where $L$ is the class of the tautological sheaf, $G$
is the class of a fiber, and
$$
V\sim 3L+mG
$$
is a smooth sufficiently general divisor in the linear
system $|3L+mG|$. Assuming that $3L+mG$ is an ample class,
we get by the Lefschetz theorem that
$$
\mathop{\rm Pic} V={\mathbb Z} K_V\oplus {\mathbb Z} F,
$$
where $F=G|_V$ is the class of a fiber. Furthermore,
$$
A^2V={\mathbb Z} s\oplus {\mathbb Z} f,
$$
where $s$ is the class of some section and $f$ is the
class of a line in a fiber.
\vspace{0.5cm}

{\bf Theorem [P4].} {\it Assume that $K^2_V$ does not
lie in the interior of the cone of effective curves
in $A^2V$. Then $V/{\mathbb P}^1$ is birationally rigid.
}
\vspace{0.5cm}

{\bf Remark.} The group of proper birational self-maps
$\mathop{\rm Bir} (V/{\mathbb P}^1)$ was described 
(generators and relations)
by Yu.I.Manin, see [M3]. It is very big.
\vspace{0.3cm}

{\bf Example.} Let ${\cal E}={\cal O}^{\oplus 4}_{{\mathbb P}^1}$,
so that ${\mathbb P}({\cal E})={\mathbb P}^3\times{\mathbb P}^1$. The variety $V$
is a divisor of bidegree $(3,m)$. It is easy to check
that for $m\geq 3$ the $K^2$-condition is satisfied,
so that by the theorem $V/{\mathbb P}^1$ is birationally rigid.

If $m=2$, then the projection
$$
p\colon V\to{\mathbb P}^3
$$
onto the first factor is of degree two. Therefore there
exists a Galois involution
$$
\tau\in \mathop{\rm Bir} V,
$$
which by construction can not be fiber-wise. Using the
method of [P4], it was proved in [Sob] that the group
$\mathop{\rm Bir} V$ is isomorphic to the {\it free} product
$$
\mathop{\rm Bir} V=\mathop{\rm Bir} F_{\eta}\mathop{*}\langle \tau\rangle
$$
and $V/{\mathbb P}^1$ is ``almost rigid'' in the following 
sense: each pencil of rational surfaces on $V$ can be
transformed into the pencil of fibers $|F|$ by means
of a birational self-map. There are no conic bundle
structures on $V$. In particular, $V$ is non-rational.
(Another example of ``almost rigidity'' see in [Gr].)

Combining the $m\geq 3$ and $m=2$ cases, we get a 
far-reaching generalization of the remarkable theorem
of late Fabio Bardelli [B], obtained by means of the
Clemens degeneration method.

The above-described results are quite precise: all 
the cases
$m\geq 2$ are embraced, while if $m=1$,  then $V$ is
clearly rational.
\vspace{1.5cm}


{\large\bf References }
\vspace{0.5cm}

\noindent
[B] Bardelli F., Polarized mixed Hodge structures: On irrationality of threefolds via degeneration. Ann. Mat. Pura et Appl. {\bf 137} (1984), 287-369.
\vspace{0.5cm} \par \noindent
[CG] Clemens H. and Griffiths P., The intermediate Jacobian of the cubic threefold. Ann. of Math. {\bf 95} (1972), 281-356. 
\vspace{0.5cm} \par \noindent
[C1] Corti A., Factoring birational maps of threefolds after Sarkisov. J. Algebraic Geom. {\bf 4} (1995), no. 2, 
223-254.
\vspace{0.5cm} \par \noindent
[C2] Corti A., Singularities of linear systems and 3-fold
birational geometry, in ``Explicit Birational Geometry
of Threefolds'', London Mathematical Society Lecture Note
Series {\bf 281} (2000), Cambridge University Press, 259-312.
\vspace{0.5cm} \par \noindent
[CM] Corti A. and Mella M., Birational geometry of terminal
quartic 3-folds. I, preprint, math.AG/0102096
\vspace{0.5cm}\par\noindent
[CPR] Corti A., Pukhlikov A. and Reid M., Fano 3-fold
hypersurfaces, in ``Explicit Birational Geometry
of Threefolds'', London Mathematical Society Lecture Note
Series {\bf 281} (2000), Cambridge University Press, 175-258.
\vspace{0.5cm} \par \noindent
[CR] Corti A. and Reid M., Foreword to
``Explicit Birational Geometry
of Threefolds'', London Mathematical Society Lecture Note
Series {\bf 281} (2000), Cambridge University Press, 1-20.
\vspace{0.5cm} \par \noindent
[F1] Fano G., Sopra alcune varieta algebriche a tre dimensioni aventi tutti i generi nulli, Atti Acc. Torino {\bf 43} (1908), 973-977.
\vspace{0.5cm} \par \noindent
[F2] Fano G., Osservazioni sopra alcune varieta non razionali aventi tutti i generi nulli, Atti Acc. Torino 
{\bf 50} (1915), 1067-1072.
\vspace{0.5cm} \par \noindent
[F3] Fano G., Nouve ricerche sulle varieta algebriche a tre dimensioni a curve-sezioni canoniche, Comm. Rent. Ac. Sci. {\bf 11} (1947), 635-720.
\vspace{0.5cm} \par \noindent
[G] Gizatullin M. Kh., Defining relations for the Cremona group of the plane. Izv. Akad. Nauk SSSR Ser. Mat. {\bf 46} (1982), no. 5, 909-970, 1134. 
\vspace{0.5cm} \par \noindent
[Gr] Grinenko M. M., Birational automorphisms of a three-dimensional double cone, Mat. Sb. {\bf 189} (1998), 
no. 7, 37-52 (English translation in Sb. Math. {\bf 189} (1998), no. 7-8, 991-1007).
\vspace{0.5cm} \par \noindent
[H] Hudson H., Cremona transformations in plane and
space, Cambridge Univ. Press, Cambridge 1927.
\vspace{0.5cm} \par \noindent
[I1] Iskovskikh V.A., Rational surfaces with a pencil
of rational curves, Mat. Sb. {\bf 74} (1967), 608-638
(English translation in Math. USSR-Sbornik {\bf 3} (1967)).
\vspace{0.5cm} \par \noindent
[I2] Iskovskikh V.A., Rational surfaces with a pencil
of rational curves and with positive square of 
canonical class, Mat. Sb. {\bf 83} (1970), 90-119
(English translation in Math. USSR-Sbornik {\bf 12} (1970)).
\vspace{0.5cm} \par \noindent
[I3] Iskovskikh V.A., Birational automorphisms of
three-dimensional algebraic varieties, J. Soviet Math.
{\bf 13} (1980), 815-868.
\vspace{0.5cm} \par \noindent
[I4] Iskovskikh V.A., A simple proof of Gizatullin's
theorem on relations in the two-dimensional Cremona
group, Trudy Mat. Inst. Steklov {\bf 183} (1990), 111-116
(English translation in Proc. Steklov Inst. Math. 1991,
no. 4).
\vspace{0.5cm} \par \noindent
[I5] Iskovskikh V. A. Factorization of birational mappings of rational surfaces from the point of view of Mori theory. Uspekhi Mat. Nauk {\bf 51} (1996), no. 4, 3-72 (English
translation in Russian Math. Surveys {\bf 51} (1996), no. 4, 585-652).
\vspace{0.5cm} \par \noindent
[IM] Iskovskikh V.A. and Manin Yu.I., Three-dimensional
quartics and counterexamples to the L\" uroth problem,
Math. USSR Sb. {\bf 86} (1971), no. 1, 140-166.
\vspace{0.5cm} \par \noindent
[IP] Iskovskikh V.A. and Pukhlikov A.V., Birational
automorphisms of multi-dimensional algebraic varieties,
J. Math. Sci. {\bf 82} (1996), 3528-3613.
\vspace{0.5cm} \par \noindent
[K] Koll{\'a}r J., et al., Flips and Abundance for
Algebraic Threefolds, Asterisque 211, 1993.
\vspace{0.5cm} \par \noindent
[M1] Manin Yu. I., Rational surfaces over perfect fields.  Publ. Math. IHES {\bf 30} (1966), 55-113.
\vspace{0.5cm} \par \noindent
[M2] Manin Yu. I. Rational surfaces over perfect fields. II.  Mat. Sb. {\bf 72} (1967), 161-192.
\vspace{0.5cm} \par \noindent
[M3] Manin Yu. I., Cubic forms. Algebra, geometry, arithmetic. Second edition. North-Holland Mathematical Library, {\bf 4.} North-Holland Publishing Co., Amsterdam, 1986.
\vspace{0.5cm} \par \noindent 
[N] Noether M., {\" U}ber Fl{\" a}chen welche Schaaren rationaler Curven besitzen, Math. Ann. {\bf 3} (1871), 161-227.
\vspace{0.5cm} \par \noindent
[P1] Pukhlikov A.V., Birational automorphisms of a three-dimensional quartic with an elementary singularity, Math. USSR Sb. {\bf 63} (1989), 457-482.
\vspace{0.5cm} \par \noindent
[P2] Pukhlikov A. V., Maximal singularities on the Fano variety $V_6^3$. Vestnik Moskov. Univ. Ser. I Mat. Mekh. 1989,  no. 2, 47-50 (English translation in Moscow
Univ. Math. Bull. {\bf 44} (1989), no. 2, 70-75).
\vspace{0.5cm} \par \noindent
[P3] Pukhlikov A.V., Essentials of the method of maximal
singularities, in ``Explicit Birational Geometry
of Threefolds'', London Mathematical Society Lecture Note
Series {\bf 281} (2000), Cambridge University Press, 73-100.
\vspace{0.5cm} \par \noindent
[P4] Pukhlikov A.V., Birational automorphisms of
three-dimensional algebraic varieties with a pencil of
del Pezzo surfaces, Izvestiya: Mathematics {\bf 62}:1 (1998), 115-155.
\vspace{0.5cm} \par \noindent
[P5] Pukhlikov A.V., Birational automorphisms of Fano
hypersurfaces, Invent. Math. {\bf 134} (1998), no. 2,
401-426. \vspace{0.5cm} \par \noindent
[P6] Pukhlikov A.V., Birationally rigid Fano double
hypersurfaces, Sbornik: Mathematics {\bf 191} (2000), No. 6,
101-126.
\vspace{0.5cm} \par \noindent
[P7] Pukhlikov A.V., Birationally rigid Fano fibrations,
Izvestiya: Mathematics {\bf 64} (2000), 131-150.
\vspace{0.5cm}\par \noindent
[P8] Pukhlikov A.V., Birationally rigid Fano complete
intersections, Crelle J. f\" ur die reine und angew. Math. {\bf 541} (2001), 55-79. \vspace{0.5cm}\par\noindent 
[P9] Pukhlikov A.V., Birationally rigid Fano
hypersurfaces with isolated singularities, Sbornik: Mathematics  {\bf 193} (2002), no. 3, 445-471. \vspace{0.5cm}
\par \noindent
[P10] Pukhlikov A.V., Birationally rigid Fano hypersurfaces,
Izvestiya: Mathematics {\bf 66} (2002), no. 6; 
math.AG/0201302.
\vspace{0.5cm}\par \noindent
[R] Reid M., Birational geometry of 3-folds according to
Sarkisov. Warwick Preprint, 1991.
\vspace{0.5cm} \par \noindent
[S1] Sarkisov V.G., Birational automorphisms of conic
bundles, Izv. Akad. Nauk SSSR, Ser. Mat. {\bf 44} (1980),
no. 4, 918-945 (English translation:
Math. USSR Izv. {\bf 17} (1981), 177-202).
\vspace{0.5cm} \par \noindent
[S2] Sarkisov V.G., On conic bundle structures, Izv. Akad.
Nauk SSSR, Ser. Mat. {\bf 46} (1982), no. 2, 371-408
(English translation: Math. USSR Izv. {\bf 20} (1982), no. 2, 354-390).
\vspace{0.5cm} \par \noindent
[S3] Sarkisov V.G., Birational maps of standard ${\mathbb Q}$-Fano fibrations, Preprint, Kurchatov Institute of
Atomic Energy, 1989.
\vspace{0.5cm} \par \noindent
[Sh] Shokurov V.V., 3-fold log flips, Izvestiya:
Mathematics {\bf 40} (1993), 93-202. \vspace{0.5cm}
\par \noindent 
[Sob] Sobolev I. V., Birational automorphisms of a class of varieties fibered into cubic surfaces. Izv. Ross. Akad. 
Nauk Ser. Mat. {\bf 66} (2002), no. 1, 203-224.
\vspace{0.5cm} \par \noindent

\end{document}